\documentclass[12pt]{amsart}

\usepackage{graphicx}
\usepackage{times}
\newtheorem{thm}{Theorem}[section] \newtheorem{cor}[thm]{Corollary}
 
\newtheorem{defin}[thm]{Definition} \newtheorem{lemma}[thm]{Lemma}
 \newtheorem{prop}[thm]{Proposition}
\newtheorem{rmk}[thm]{Remark}

\newcommand{\mcD}{\mbox{$\mathcal D$}}

\newcommand{\mcB}{\mbox{$\mathcal B$}}

\newcommand{\mcH}{\mbox{$\mathcal H$}}

\newcommand{\bdd}{\mbox{$\partial$}}

\begin{document}



\title{Thin Position for knots in a 3-manifold}

\author{Maggy Tomova} \address{\hskip-\parindent Maggy Tomova\\
  Mathematics Department \\
University of Iowa \\
 Iowa city, IA 52240, USA} \email{mtomova@math.uiowa.edu}

\date{\today}

\begin{abstract}
We extend the notion of thin multiple Heegaard splittings of a link 
in a 3-manifold to take into consideration not only compressing disks 
but also cut-disks for the Heegaard 
surfaces. We prove that if $H$ is a c-strongly compressible bridge surface 
for a link $K$ contained in a closed orientable irreducible 3-manifold $M$ then one 
of the following is satisfied:

\begin{itemize}
    
    \item $H$ is stabilized
    \item $H$ is meridionally stabilized
    \item $H$ is perturbed
    \item a component of $K$ is removable
    \item $M$ contains an essential meridional surface. 
 
 \end{itemize}   
  
\end{abstract}

\maketitle
\section{Introduction}

The notion of thin position for a closed orientable 3-manifold $M$ was introduced by 
Scharlemann and Thompson in \cite{ST1}. The idea is to build the 
3-manifold by starting with a set of 0-handles, then alternate 
between attaching collections of 1-handles and 2-handles keeping the 
boundary at the intermediate steps as simple as possible and finally 
add 3-handles. Such a decomposition of a manifold is called a 
generalized Heegaard splitting. The classical Heegaard splitting 
where all 1-handles are attached at the same time followed by all 
2-handles is an example of a generalized Heegaard splitting. 
Casson and Gordon \cite{CG} show that if $A\cup_P B$ is a weakly reducible 
Heegaard splitting for $M$ (i.e. there are meridional disks for $A$ 
and $B$ with disjoint boundaries), then either $A\cup_P B$ is 
reducible or $M$ contains an essential surface. Scharlemann and 
Thompson \cite{ST1} show that such surfaces arise naturally 
when a Heegaard splitting in put in thin position.

Suppose a closed orientable 3-manifold $M=A\cup_P B$ contains a link 
$K$, then we can 
isotope $K$ so that it intersects each handlebody in boundary parallel arcs. 
In this case we say that $P$ is a bridge surface for $K$ or that $P$ is a 
Heegaard surface for the pair $(M,K)$. 
The idea was first introduced by Schubert in the case that $M=S^3$ 
and $P=S^2$ and was extended by Morimoto and Sakuma for other 
3-manifolds. In \cite{HS1} Hayashi and Shimokawa considered multiple 
Heegaard splittings for $(M,K)$ using the idea of changing the order 
in which the 1-handles and the 2-handles are attached. They 
generalized the result of \cite{ST1} in this context, i.e. 
they showed that 
if $P$ is a strongly compressible bridge surface for $K$, then  
either $A\cup_P B$ is stabilized or cancellable or $M-\eta(K)$ contains an 
essential meridional surface.  

In this paper we will generalize this important result one step 
further by weakening the hypothesis. Suppose $M$ is a compact
orientable manifold and $F\subset M$ is a properly embedded surface transverse 
to a 
1-submanifold $T\subset M$. In some contexts it is necessary to consider not only compressing 
disks for $F$ but also cut-disks, that is, disks whose boundary is 
essential on $F-T$ and that intersect $T$ in exactly one point, as for 
example in \cite{BS}, \cite{STo2} and \cite{T2}. A bridge surface $P$ for 
a link $K$ is 
c-strongly compressible if there is a pair of disjoint cut or compressing 
disks for $P_K$ on 
opposite sides of $P$. In particular every strongly compressible 
bridge surface is c-strongly compressible. We will show that if a 
bridge surface $P$ for $K$
is c-strongly compressible then either it can be simplified in one of four 
geometrically obvious ways or $(M,K)$ contains an essential meridional surface.

\section{Definitions and preliminaries}

Let $M$ be a compact orientable irreducible 3-manifold and let $T$ be a 1-manifold properly 
embedded in $M$. A regular neighborhood of $T$ will be denoted 
$\eta(T)$. If $X$ is any subset of $M$ we will use $X_T$ to denote $X-T$.
We will assume that any sphere in $M$ intersects $T$ in an even number 
of points. As all the results we will develop are used in the 
context when $T$ only has closed components, this is a natural 
assumption. If $K$ is a link in $M$, then any sphere in $M$ intersects $K$ in an even number of 
points, since the ball in bounds in $M$ contains no endpoints of $K$.

Suppose $F$ is a properly embedded surface in $M$. An {\em essential curve} on 
$F_T$ 
is a curve that doesn't bound a disk on $F_T$ and it is not parallel to 
a puncture of $F_T$. A {\em compressing disk} $D$ for $F_T$ is an embedded disk in 
$M_T$ so that 
$F \cap D =\bdd D$ is an essential curve on $F_T$. A {\em cut-disk} is a 
disk $D^c \subset M$ such that $D^c \cap F=\bdd D^c$ is an essential 
curve on $F_T$ and $|D \cap T|=1$. A {\em c-disk} is a cut or a 
compressing disk. $F$ will be called {\em incompressible} if it has no 
compressing disks and {\em c-incompressible} if it has no c-disks. $F$ will be called {\em essential} if it does not have compressing disks 
(it may have cut disks), it is not boundary parallel in $M-\eta(T)$ and it 
is not a sphere that bounds a ball disjoint from $T$.

Suppose $C$ is a compression body ($\bdd_-C$ may have some sphere 
components). A set of arcs $t_i \subset C$ is {\em trivial} if there is a 
homeomorphism after which each arc is either vertical, ie, 
$t_i=(point) \times I \subset \bdd_-C \times I$ or there is an 
embedded disk $D_i$ such that $\bdd D_i=t_i \cup \alpha_i$ where 
$\alpha_i \subset \bdd_+C$. In the second case we say that $t_i$ is 
$\bdd_+$-parallel and the disk $D_i$ is a bridge disk. If $C$ is a 
handlebody, then all trivial arcs are $\bdd_+$-parallel and are called 
{\em bridges}.  If $T$ is a 1-manifold properly embedded in a 
compression body $C$ so that $T$ is a collection of trivial arcs then 
we will denote the pair by $(C, T)$.

Let $(C,T)$ be a pair of a compression body and a 1-manifold and 
let $\mcD$ be
the disjoint union of compressing disks for $\bdd_+C$ together with 
one bridge disk for each $\bdd_+$-parallel arc. If $\mcD$ cuts $(C,T)$ 
into a manifold homeomorphic to $(\bdd_- C\times I,$ {\em vertical 
arcs}$)$ together with some 3-balls, then $\mcD$ is called a {\em complete disk system} 
for $(C,T)$. The presence of such a complete disk system can 
be taken as the definition of $(C,T)$.

Let $M$ be a 3-manifold, let $A \cup_P B$ be a Heegaard splitting 
(ie $A$ and $B$ are compression bodies) for 
$M$ and let $T$ be a 1-manifold in $M$. We say that $T$ is in bridge position with 
respect to $P$ if $A$ and $B$ intersect $T$ only in trivial arcs. In 
this case we say that $P$ is a bridge surface for $T$ or that $P$ as a Heegaard surface for the pair $(M,T)$.

Suppose $M=A \cup_ P B$ and $T$ is in bridge position with respect 
to $P$. The Heegaard splitting  is {\em c-strongly irreducible} if any 
pair of c-disks on opposite sides of $P_T$ intersect, in this case the 
bridge surface $P_T$ 
is {\em c-weakly incompressible}. If there are c-disks $D_A 
\subset A$ and $D_B \subset B$ such that $D_A \cap D_B =\emptyset$, the
Heegaard splitting is {\em c-weakly reducible} and the bridge surface 
$P_T$ is 
{\em c-strongly compressible}.

Following \cite{HS1}, the bridge surface $P_T$ will be called {\em stabilized} if there is a pair of 
compressing disks on opposite sides of $P_T$ that intersect in a single 
point. The bridge surface is {\em meridionally stabilized} if 
there is a cut disk and a compressing disk on opposite sides of $P_T$ 
that intersect in a single point. Finally the bridge surface is called 
{\em cancellable} if there is a pair of 
canceling disks $D_i$ for bridges $t_i$ on opposite sides of $P$ 
such that $ \emptyset \neq (\bdd D_1 \cap \bdd D_2) \subset (Q \cap K)$. 
If $|\bdd D_1 \cap \bdd D_2|=1$ we will call the bridge surface
{\em perturbed}. In \cite{STo3} the authors show that is $M=A \cup_P 
B$ is stabilized, meridionally stabilized or perturbed, then there is 
a simpler bridge surface $P'$ for $T$ such that $P$ can be obtained 
from $P'$ by one of three obvious geometric operations. 

If the 
bridge surface $P$ for $T$ is cancellable with 
canceling disks $D_1$ and $D_2$ such that $|\bdd D_1 \cap \bdd D_2|=2$ 
then using this pair of disks some closed component $t$ of $T$ can be isotoped to 
lie in $P$. If this component can be isotoped to lie in the core of one of the 
compression bodies, $A$ say, and is disjoint from all other bridge 
disks in $A$ then $A - \eta(t)$ is also a compression 
body and the 1-manifold $T-t$ intersects it in a collection 
of trivial arcs. Thus $(A-\eta(t)) \cup_P B$ is 
Heegaard splitting for $(M -\eta(t))$ and $P$ is a bridge surface for 
$T-t$. In this case we will say that $T$ has 
a {\em removable component}. A detailed discussion of links 
with removable components is given in \cite{STo3}.

In the absence of a knot, it follows by a theorem of Waldhausen that a 
Heegaard splitting of an irreducible manifold is stabilized if and only if there is a sphere that 
intersects the Heegaard surface in a single essential curve (i.e the 
Heegaard splitting is reducible), unless 
the Heegaard splitting is the standard genus 1 Heegaard splitting of 
$S^3$. We will say that a bridge surface for $T$ is {\em c-reducible} 
if there is a sphere or a twice punctured sphere in $M$ that intersects 
the bridge surface in a 
single essential closed curve. Then one direction of Waldhausen's result 
easily generalizes to 
bridge surfaces as the next theorem shows.

\begin{thm} \label{thm:reducible}
    Suppose $P$ is a bridge surface for a 1-manifold $T$ properly 
    embedded in a compact, orientable 3-manifold $M$ where $P$ is not 
    the standard genus 1 Heegaard splitting for $S^3$. If $P$ is stabilized, 
    perturbed or meridionally stabilized then there exists a 
    sphere $S$, possibly punctured by $T$ twice, which intersects 
    $P$ in a single essential curve $\alpha$ and neither component of 
    $S - \alpha$ is parallel to $P$. 
    
\end{thm}

\begin{figure}[tbh]
\centering
\includegraphics[scale=0.5]{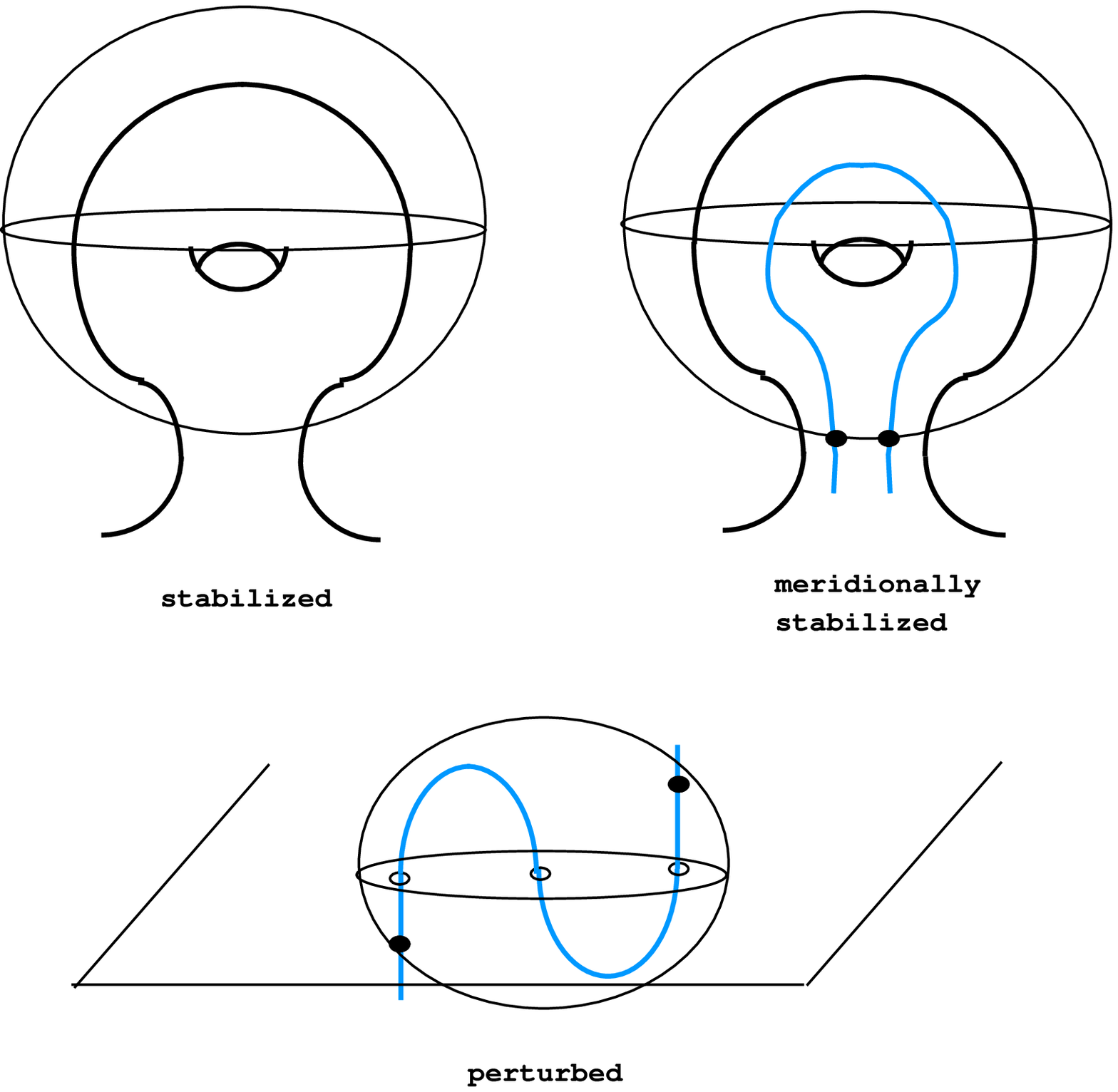}
\caption{} \label{fig:stabilizations}
\end{figure}

\begin{proof}
    If $P$ is stabilized let $S$ be the boundary of a regular neighborhood of the union of 
    the pair of stabilizing disks, Figure \ref{fig:stabilizations}. In 
    this case $S$ is a sphere disjoint from $T$. 
    If $P$ is meridionally stabilized, let $S$ be the boundary of a  
    regular neighborhood of the union of the cut and compressing 
    disks.  In this case $S$ is a twice punctured sphere with both punctures on the 
    same side of $S \cap P$. Finally if $P$ is perturbed with canceling disks $E_1$ and $E_2$, 
    let $S$ be the boundary of a regular neighborhood of $E_1 \cup E_2$. Then $S$ is a twice punctured sphere 
    and the punctures are separated by $S \cap P$.
    
     %
     %
 \end{proof}   
 
 \section{C-compression bodies and their properties}
  
 We will need to generalize the notion of a compression body containing 
 trivial arcs as follows.

\begin{defin}
    A c-compression body $(C, T)^c$ is a pair of a compression body 
    $C$ and a 1-manifold $T$ such that there is a collection of 
    disjoint bridge disks and c-disks $\mcD^c$ so that $\mcD^c$ cuts 
    $(C,T)^c$ into a 3-manifold homeomorphic to $(\bdd_-C 
 \times I,$ vertical arcs$)$ together with some 3-balls. In 
 this case $\mcD^c$ is called a complete c-disk system.
    
 \end{defin}

 One way to construct a compression body is to take 
 a product neighborhood $F \times 
 I$ of a closed, possibly disconnected, surface 
 $F$ so that any arc of $T \cap (F \times I)$ can either be isotoped to 
 be vertical with 
 respect to the product structure or is parallel to an arc in $F \times 
 0$ and then attach a collection of pairwise disjoint 2-handles $\Delta$ to
 $F \times 1$. If we allow some of the 2-handles in $\Delta$ to contain 
 an arc $t \subset T$ as their cocore, the resulting 3-manifold is a 
 c-compression body. The complete c-disk system described in the 
 definition above consists of all 
 bridge disks together with the cores of the 2-handles. We will use 
 this construction as an alternative definition of a c-compression body.  

\begin{rmk}
 {\em Recall that a spine of a compression body $C$ is the union of $\bdd_- 
 C$ together with a 1-dimensional graph $\Gamma$ such that $C$ 
 retracts to $\bdd_-C \cup \Gamma$. An equivalent definition of a c-compression body 
 is that  $(C,T)^c$ is a compression body $C$ together with 
 a 1-manifold $T$ and there exists a spine $\Sigma$ for $C$ such that 
all arcs of $T$ that are not trivial in $C$ can be simultaneously  
isotoped to lie on $\Sigma$ and be pairwise disjoint. We will however not use 
this definition here.}
 \end{rmk}

 \begin{prop}\label{prop:c-compcanbecomp}
 
  Let $(C,T)^c$ be a c-compression body. 
  Then $(C,T)^c$ is a compression body if and only if there is no arc $t \subset T$ such that $\bdd t \subset  \bdd_-C$. 
  In particular if $\bdd_-C = \emptyset$, then $C$ is a handlebody.
  \end{prop}

  \begin{proof}
 Consider the construction above and note that before the two handles 
 are added no arc of $T$ has both of its 
 endpoints on $F \times 1$. If some 2-handle $D$ attached 
 to $F \times 1$ contains an arc $t \subset T$ as its core, this arc 
 will have both of its endpoints on $\bdd_- C$. Thus $C$ is a 
 compression body if and only if no 2-handle contains such an arc.

  \end{proof}    

  \begin{lemma} \label{lem:surfaceincbody}
	 Let $(C,T)^c$ be a c-compression body and let $F$ be a c-incompressible, $\bdd$-incompressible
      properly embedded surface transverse to $T$.  Then there is a
      complete c-disk system $\mcD^c$ of $(C,T)^c$ such that $\mcD^c \cap F$
      consists of two types of arcs
      \begin{itemize}
	  \item An intersection arc $\alpha$ between a bridge disk
	  in $\mcD^c$ and a twice punctured sphere component of $F$ with both endpoints of
	  $\alpha$ lying on $T$.  
	  
	  \item An intersection arc $\beta$
	  between a bridge disk in $\mcD^c$ and a once-punctured disk
	  component of $F$ with one endpoint of $\beta$ lying on $T$ and
	  the other lying on $\bdd_+ C$.
  \end{itemize}
      
  \end{lemma}    

  \begin{proof}
   The argument is similar to the the proof of Lemma 2.2 in \cite{HS1} so we only 
   give an outline here. Let $\mcD^c$ be a complete c-disk system for 
   $(C,T)^c$ chosen to minimize $|\mcD^c \cap F|$. Using the fact that 
   $F_T$ is c-incompressible, we may assume that $\mcD^c \cap F$ does not 
   contain any simple closed curves. If $\alpha \subset \mcD^c \cap F$ is 
   an arc with both of its endpoints on $\bdd C$, then an outermost 
   such arc either gives a $\bdd$-compression for $F$ contrary to the 
   hypothesis or can be removed by an outermost arc argument 
   contradicting the minimality of $|\mcD^c \cap F|$. Note that if 
   $\alpha$ lies on some cut-disk $D^c$, we can still choose the arc so that 
   the disk it cuts from $D^c$ does not contain a puncture.  This establishes 
   that $F$ is disjoint from all c-disks in $\mcD^c$. 
   
   Suppose $\alpha$ is an 
   arc of intersection between a bridge disk $E$ for $T$ and a 
   component $F'$ of $F$. Assume that $\alpha$ is an outermost such arc 
   and let $E' \subset E$ be the subdisk $\alpha$ bounds on $E$. By 
   the above argument at least one endpoint of $\alpha$ must lie on 
   $T$. If 
   both endpoints of $\alpha$ lie on $T$, the boundary of a regular 
   neighborhood of $E'$ gives a compressing disk for $F$ contrary 
   to the hypothesis unless $F'$ is a twice punctured sphere.
   If $\alpha$ has one endpoint on $T$ and one 
   endpoint on $\bdd C$, a regular neighborhood of $E'$ is a 
   $\bdd$-compressing disk for $F$ unless $F'$ is a once punctured disk.  
   
  \end{proof}

  \begin{cor} \label{cor:negativeboundaryincomp}
      If $(C,T)^c$ is a c-compression body, then $\bdd_- C$ is 
      incompressible.
      
  \end{cor}    
  
  \begin{proof}
      Suppose $D$ is a compressing disk for some component of $\bdd_- C$. By 
      Lemma \ref{lem:surfaceincbody}, there exists a complete c-disk 
      system $\mcD^c$ for $(C,T)^c$ such that $D \cap \mcD^c = \emptyset$. 
      But this implies that $D$ is a $\bdd$-reducing disk for the 
      manifold $(F \times I, $ {\em vertical arcs}$)$, a 
      contradiction.
      
      \end{proof}

  If $M$ is a 3-manifold we will denote by $\tilde M$ the manifold 
  obtained from $M$ by filling any sphere boundary components of $M$ 
  with 3-balls. 

  \begin{lemma} [Lemma 2.4 \cite{HS1}]\label{lem:surfaceinproduct}
      If $F$ is an incompressible, $\bdd$-incompressible surface in a 
	 compression body $(C, T)$, then $F$ is a collection of the following 
	 kinds of components:
	 
	\begin{itemize}
	    \item Spheres intersecting $T$ in 0 or 2 points,
	    \item Disks intersecting $T$ in 0 or 1 points,
	    \item Vertical annuli disjoint from $T$,
	    \item Closed surfaces parallel to a component of $\bdd_- \tilde C$.
	    \end{itemize}
      \end{lemma}

  \begin{cor} \label{cor:components}
      If $F$ is a c-incompressible, $\bdd$-incompressible surface in a 
      c-compression body $(C,T)^c$, then $F$ is a collection of the following 
      kinds of components:
      
     \begin{itemize}
	 \item Spheres intersecting $T$ in 0 or 2 points,
	 \item Disks intersecting $T$ in 0 or 1 points,
	 \item Vertical annuli disjoint from $T$,
	 \item Closed surfaces parallel to a component of $\bdd_- \tilde C$.
	 \end{itemize}
      
   \end{cor}   

   \begin{proof}
       Delete all component of the first two types and let $F'$ be the 
       new surface. By Lemma \ref{lem:surfaceincbody}, there exists a complete 
       c-disk system $\mcD$ for $(C,T)^c$ such that $\mcD \cap F' = \emptyset$. 
       Thus each component of $F'$ is contained in a compression body 
       with trivial arcs (in fact in a trivial compression body but we don't need this 
       fact). The result follows by Lemma \ref{lem:surfaceinproduct}.

   \end{proof}

 \section {C-thin position for a pair 3-manifold, 1-manifold}
 
The following definition was first introduced in \cite{HS1}

\begin{defin}If $T$ is a 1-manifold properly embedded in a compact 
3-manifold $M$, we say that the disjoint union of surfaces $\mcH$ is 
a multiple Heegaard splitting of $(M,T)$ if 
\begin{enumerate}
\item The closures of all components of $M-\mcH$ are compression 
bodies $(C_1, C_1\cap T),...,(C_n, C_n \cap T)$,
\item for $i=1,...,n$,  $\bdd_+ C_i$ is attached to some $\bdd_+C_j$  
where $i\neq 
j$,
\item a component of $\bdd_- C_i$ is attached to some component of 
$\bdd_-C_j$ (possibly $i=j$).

\end{enumerate}
A component $H$ of $\mcH$ is said to be positive if $H=\bdd_+C_i$ for 
some $i$ and negative if $H=\bdd_-C_j$ for some $j$. The unions of 
all positive and all negative components of $\mcH$ are denoted 
$\mcH_+$ and $\mcH_-$ respectively.

\end{defin}

Note that if $\mcH$ has a single surface component $P$, then $P$ is a bridge surface for $T$.

Using c-compression bodies instead of compression bodies, we 
generalize this definition as follows.

\begin{defin}If $T$ is a 1-manifold properly embedded in a compact 
3-manifold $M$, we say that the disjoint union of surfaces $\mcH$ is 
a multiple c-Heegaard splitting of $(M,T)$ if 
\begin{enumerate}
\item The closures of all components of $M-\mcH$ are c-compression 
bodies $(C_1, C_1\cap T)^c,...,(C_n, C_n \cap T)^c$,
\item for $i=1,...,n$, $\bdd_+ C_i$ is attached to some $\bdd_+C_j$ 
where $i\neq 
j$,
\item a component of $\bdd_- C_i$ is attached to some component of 
$\bdd_-C_j$ (possibly $i=j$)
\end{enumerate}

\end{defin}

As in \cite{ST1} and \cite{HS1} we will associate to a multiple 
c-Heegaard splitting a measure of its complexity. The following 
notion of complexity of a surface is different from the 
one used in \cite{HS1}.

\begin{defin}Let $S$ be a closed connected surfaces 
embedded in $M$ transverse to a properly embedded 1-manifold $T 
\subset M$. The complexity of $S$ is the 
ordered pair $c(S)=(2-\chi(S_T), g(S))$. If $S$ is not connected, 
$c(S)$ is the multi-set of ordered pairs corresponding to each of the 
components of $S$.
\end{defin}

As in 
\cite{ST1} the complexities 
of two possibly not connected surfaces are compared by first arranging the 
ordered pairs in each multi-set in non-increasing order and then 
comparing the two multi-sets lexicographically where the ordered pairs 
are also compared lexicographically. 

\begin{lemma} \label{lem:compreduces}
Suppose $S_T$ is meridional surface in $(M,T)$ of non-positive euler 
characteristic. If $S_T'$ is a component of the surface obtained from $S_T$ by compressing along 
a c-disk, then $c(S_T)>c(S'_T)$.
\end{lemma}
 
\begin{proof} Without loss of generality we may assume that $S_T$ is 
connected.

{\em Case 1:} Let $\tilde S_T$ be a possibly disconnected surface obtained from $S_T$ via a compression along a 
disk $D$. In this case $\chi(S_T) < \chi(\tilde S_T)$ as $\chi(D)=1$ 
so the result follows immediately if $\tilde S_T$ is connected. If 
$\tilde S_T$ consists of two components then by the definition of compressing 
disk, we may assume that neither component is a sphere and thus both
components of $\tilde S_T$ have non-positive Euler characteristic. By 
the additivity of Euler characteristic it 
follows that if $S'_T$ is a component of $\tilde S_T$, then $\chi( \tilde 
S_T)\leq \chi(S'_T)$ so $2-\chi(S_T)>2-\chi(S'_T)$ as desired.

{\em Case 2:} Suppose $\tilde S_T$ is obtained from $S_T$ via a compression along a 
cut-disk $D^c$. If $D^c$ is separating, then each of the two components of $\tilde S_T$ has at least 
one puncture and if a component is a sphere, then it must have at 
least 3 punctures, ie each component of $\tilde S_T$ has a strictly 
negative Euler characteristic. By the additivity of Euler 
characteristic, we conclude that for each component $S'_T$ of $\tilde 
S_T$, $\chi(S'_T)< \chi(\tilde S_T) = \chi(S_T)$ and so the first component 
of the complexity tuple is decreased.

If the cut disk is not separating the cut-compression does not 
affect the first term in the complexity tuple as $\chi(D^c)=0$. Note 
that $\bdd D^c$ must be essential in the non-punctured surface $S$ so we can consider $D^c$ as a 
compressing disk for $S$ in $M$. Then $g(\tilde S) < g(S)$ so in this 
case the 
second component of the complexity tuple is decreased.

\end{proof}

The {\em width} of a c-Heegaard splitting is the multiset of pairs 
$w(\mcH)= c(\mcH_+)$.  In \cite{HS1} a multiple Heegaard 
splitting is called {\em thin} if it is of minimum width amongst all 
possible multiple Heegaard splittings for the pair $(M,T)$. Similarly we will 
call a c-Heegaard splitting {\em c-thin} if it is of minimal width 
amongst all c-Heegaard splittings for $(M,T)$.

\section{Thinning using pairs of disjoint c-disks}

 \begin{lemma} \label{lem:thinning}
 Suppose $M$ is a compact orientable irreducible manifold and $T$ is a
 properly embedded 1-submanifold.  If $P$ is a c-Heegaard splitting for 
 $(M,T)$ which 
 is c-weakly reducible, then there exists
 a multiple c-Heegaard splitting $\mcH'$ so that $w(\mcH') <w(P)$. 
 
 Moreover if $M$ is closed then either
 \begin{itemize}
     \item There is a component of $\mcH'_T$ that is neither an 
     inessential sphere nor boundary parallel in $M_T$, or
     \item $P$ is stabilized, meridionally stabilized or 
     perturbed, or a closed component of $T$ is removable.
     
 \end{itemize}    
 \end{lemma}

The first part of the proof of this lemma is very similar to the proof of Lemma 2.3 in
\cite{HS1} and uses the idea of {\em untelescoping}.  However, in Lemma
2.3 the authors only allow untelescoping using disks while we also
allow untelescoping using cut-disks.

 \begin{proof}

Let $(A, A\cap T)^c$ and $(B, B\cap T)^c$ be the two c-compression 
bodies that $P$ cuts $(M, T)$ into. Consider a maximal collection of c-disks $\mcD^*_A\subset A_T$ and
$\mcD^*_B\subset B_T$ such that $\bdd \mcD^*_A\cap \bdd \mcD^*_B =
\emptyset$.  Let $A'_T =cl(A_T-N(\mcD^*_A))$ and $B'_T
=cl(B_T-N(\mcD^*_B))$ where $N(\mcD*)$ is a collar of $\mcD^*$.  
Then $A'_T$ and $B'_T$ are each the disjoint
union of c-compression bodies.  Take a small collar $N(\bdd_+A'_T)$ of
$\bdd_+A'_T$ and $N(\bdd_+B'_T)$ of $\bdd_+B'_T$.  Let
$C^1_T=cl(A'_T-N(\bdd_+A'_T))$, $C^2_T=N(\bdd_+A'_T)\cup N(\mcD^*_B)$,
$C^3_T=N(\bdd_+B'_T)\cup N(\mcD^*_A)$ and
$C^4_T=cl(B'_T-N(\bdd_+B'_T))$.  This is a new multiple c-Heegaard
splitting of $(M,T)$ with positive surfaces $\bdd_+ C_1$ and $\bdd_+
C_2$ that can be obtained from $P$ by c-compressing along $\mcD^*_A$
and $\mcD^*_B$ respectively and a negative surface $\bdd_- C_2=
\bdd_-C_3$ obtained from $P$ by compressing along both sets of
c-disks.  By Lemma \ref{lem:compreduces} it follows that $w(\mcH')
<w(P)$.

To show the second part of the lemma, suppose $A\cup_P B$ is not stabilized, 
meridionally stabilized or perturbed and no component of $T$ is 
removable and, by way of contradiction, suppose that every component of 
$\bdd_-C_2$ is a 
sphere bounding a ball that intersects $T$
in at most one trivial arc or a torus that bounds a solid
torus such that $t \subset T$ is a core curve of it.

Let $\Lambda_A$ and $\Lambda_B$ be the arcs that are the cocores of
the collections of c-disks $\mcD^*_A$ and $\mcD^*_B$ respectively.  If
$D^c$ is a cut-disk, we take $\lambda \subset T$ as its cocore.  Let
$\Lambda=\Lambda_A\cup \Lambda_B$ and note that $P$ can be recovered
from $\bdd_-C_3$ by surgery along $\Lambda$.  As $P$ is connected,
at least one component of $\bdd_-C_3$ must be adjacent to both
$\Lambda_A$ and $\Lambda_B$, call this component $F$. Unless $F_T$ is 
is an inessential sphere or boundary parallel in $M_T$ we are done. 
If $F_T$ is an inessential sphere, then by Waldhausen's result the 
original Heegaard splitting is stabilized. As $\bdd M=\emptyset$ by 
hypothesis, the remaining possibility is that $F_T$ 
is parallel in $M_T$ to part of $T$; since $F_T$ is connected it is 
either a torus bounding a solid torus with a component of $T$ as its 
core or $F_T$ is an annulus, parallel to a subarc of $T$. That is $F$ 
bounds a ball which $T$ intersects in a trivial arc.

Let $\mcB$ be the ball or solid torus $F$ bounds. 
We will assume that $\mcB$ lies on the side of $F$ that is adjacent to 
$\Lambda_A$ and that $F$ is innermost in the sense that $\mcB \cap 
\Lambda_B = \emptyset$.

Let $H = \bdd_-C_3 \cap \mcB$ and let $A'$ be the c-compression body
obtained by adding the 1-handles corresponding to the arcs $\Lambda_A \cap 
\mcB$ to a collar of $H $. (Some of these 1-handles might have subarcs 
of $T$ as their core).  Let
$B'=\mcB - A'$.  Notice that $B'$ can be obtained from $B$ by
c-compressing along all c-disks whose cocores are adjacent to $F$ and
thus $B'$ is a c-compression body.  In fact $\bdd_-B' =\emptyset$ so
$B'$ is a handlebody, let $H' = \bdd B'$.  Then $A' \cup_{H'} B'$ is a c-Heegaard splitting for
$\mcB$ decomposing in into a c-compression body $A'$ and a 
handlebody $B'$. There are two cases to consider: $\mcB$ being a 
ball intersecting $T$ is a trivial arc and $\mcB$ being a torus. We will consider 
each case separately and prove that $A' \cup_{H'} B'$ is actually a 
Heegaard splitting for $\mcB$ (i.e. $A'$ is a compression body) so we can apply previously known results.

\begin{figure}[tbh]
\centering
\includegraphics[scale=0.5]{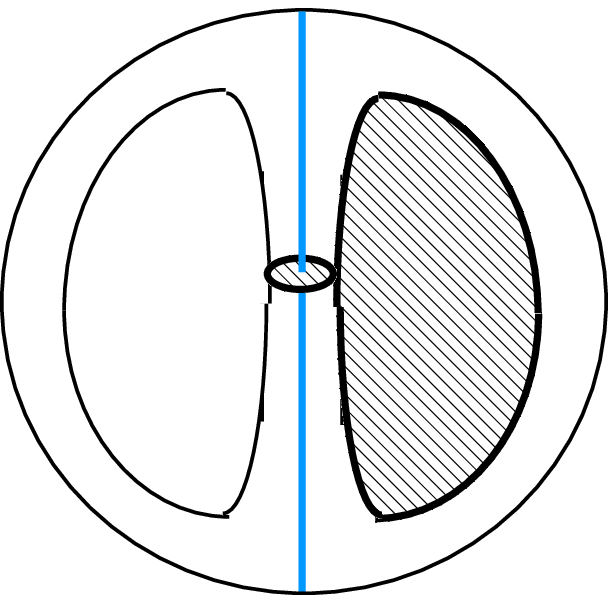}
\caption{} \label{fig:merstab}
\end{figure}

{\bf Case 1:} If $\mcB$ is a ball and $\mcB \cap T=t$ is a trivial arc, there are
three sub-cases to consider. If $t \cap H' \neq \emptyset$ then the construction 
above gives a nontrivial Heegaard splitting for the pair
$(\mcB,t)$; $A'$ is a compression body by Proposition \ref{prop:c-compcanbecomp} 
as $\bdd_-A'$ adjacent to 
two subarcs of $t$ both of which have their second endpoint on 
$\bdd_+A'=H'$.  By Lemma
2.1 of \cite{HS2}, $H'$ is either stabilized or perturbed (in this context 
if $H'$ is cancellable, it must be perturbed as $t$ is not
closed) so the same is
true for $P$.

If $t \subset A'$ and $t=\Lambda \cap \mcB$ (in particular $H =F$), Figure \ref{fig:merstab} shows a pair of c-disks
demonstrating that $P$ is meridionally stabilized.  

If $t \subset A'$ and $t\neq \Lambda \cap \mcB$, consider the solid torus
$V=\mcB-\eta(t)$.  Let $A''$ be the c-compression body obtained by
1-surgery on $H$ along the arcs $\Lambda \cap 
V$. As $t \cap V = \emptyset$, 
$A''$ is in fact a compression body. Note that
$V-A''=B'$ as $B' \cap t=\emptyset$.  Thus $A'' \cup B'$ is a
non-trivial Heegaard splitting for the solid torus $V$.  By
\cite{ST2} it must be stabilized and thus so is $P$.

{\bf Case 2:} Suppose $F$ bounds a solid torus $\mcB$, which is a regular 
neighborhood of closed component $t$ of $T$. As $\bdd_-A \cap t = 
\emptyset$, $A' \cup_{H'} B'$ is a Heegaard splitting for $(V,t)$.  By \cite{HS3}
it is cancellable or stabilized. This proves the theorem at hand unless $H'$ is cancellable but not 
perturbed so assume this is the case. In particular this implies that $H' \cap T =2$. In this case \cite{HS3} shows 
that if $g(H')
\geq 2$ then $H'$ is stabilized.  Thus it remains to consider the case 
when $H'$ is a torus intersecting $t$ in two points. In this case $H$ 
must be the union of $F$ and a sphere $S$ intersecting $t$ in two 
points and $\Lambda \cap \mcB$ is a single possibly knotted arc with one endpoint on 
$F$ and the other on $S$. As $t$ is cancellable, we can use the 
canceling disk in $A'$ to isotope $t$ across $H'$ so it lies entirely 
in $B'$. After this isotopy it is clear that $F$ and $H'$ cobound a 
product region. As $F$ is the boundary of a regular neighborhood of 
$t$, it follows that $t$ is isotopic to the core loop of the solid torus 
$B'$ ie, $B'-\eta(t)$ is a trivial 
compression body. $B$ can be recovered from $B'$ by 1-surgery so $B-\eta(t)$ is also a compression body. 
Thus after an isotopy of $t$ along the pair of canceling disks, $P$ is 
 a Heegaard splitting for $(M-\eta(t), T-t)$ so $t$ is a removable 
 component of $T$.

 \end{proof}

\section{Intersection between a boundary reducing disk and a bridge 
surface}

 As in Jaco \cite{J} a weak hierarchy for a compact orientable 
 2-manifold $F$ is a sequence of pairs $(F_0, \alpha_0),\ldots,(F_n, 
 \alpha_n)$ where $F_0=F$, $\alpha_i$ is an essential curve on $F_i$ 
 and $F_{i+1}$ is obtained from $F_i$ by cutting $F_i$ along $\alpha_i$. 
 The final surface in the hierarchy, $F_{n+1}$, satisfies the 
 following:
 
 \begin{enumerate}
     \item Each component of $F_{n+1}$ is a disc or an annulus at 
     least one boundary component of which is a component of $\bdd F$.
     \item Each non-annulus component of $F$ has at least one boundary 
     component which survives in $\bdd F_{n+1}$. 
    \end{enumerate} 
 
    The following lemma was first proven by Jaco and then extended 
    in \cite{HS1}, Lemma 3.1. 
    
    \begin{lemma} \label{lem:boundarynumber}
    Let $F$ be a connected planar surface with $b \geq 2$ boundary components. 
    Let $(F_0, \alpha_0),\ldots,(F_n, \alpha_n)$  be a weak hierarchy 
    with each $\alpha_i$ an arc. If $d$ is the number of boundary 
    components of $F_{n+1}$ then,
    \begin{itemize}
	\item If $F_{n+1}$ does not have annulus components then $d \leq b-1$ 
	\item If $F_{n+1}$ contains an annulus component, then $d \leq b$. 
	If $d=b$ and $b \geq 3$, then $F_{n+1}$ contains a disc component.
	\end{itemize}

\end{lemma}

 \begin{thm} \label{thm:one}
     Suppose $M$ is a compact orientable irreducible manifold and $T$ 
     is a properly embedded 
    1-manifold in $M$. Let $A \cup_P B$ be a c-Heegaard splitting for 
    $(M,T)$. If $D$ is a boundary reducing disk 
    for $M$ then there 
    exists such disk $D'$ so that $D'$ intersects $P_T$ in a unique 
   essential simple closed curve. 
  \end{thm}   
  
  \begin{proof}
  
       Let $D$ be a reducing disk for $\bdd M$ chosen 
      amongst all such disks so that $D \cap 
      P$ is minimal. By Corollary \ref{cor:negativeboundaryincomp}, $D \cap P \neq \emptyset$. 
      Let $D_A=D \cap A$ and $D_B=D\cap B$. 
      
      Suppose some component of $D_A$ is c-compressible in $A$ with $E$ 
      the c-compressing disk. Let $\gamma =\bdd E$ and let $D_{\gamma}$ 
      be the disk $\gamma$ bounds on $D$. Note that the sphere 
      $D_{\gamma} \cup E$ must be punctured by $T$ either 0 or two 
      times thus  $E$ must be a non-punctured disk.   
      Let $D' = (D - D_{\gamma}) \cup E$.  $D'$ is also a reducing disk for $\bdd 
      M$ as $\bdd D' =\bdd D$ and $D' \cap T= \emptyset$.  As $\bdd E$ is 
      essential on $D_A$, $D_\gamma$ cannot lie entirely in $A$ so 
      $|D_\gamma \cap P| >  |E \cap P|$ and thus $|D' \cap P| < |D \cap 
      P|$ contradicting the choice of $D$. Similarly $D_B$ is 
      c-incompressible in $B$.
      
       Suppose that $E$ is a $\bdd$-compressing disk for $D_A$ and $E$ 
      is adjacent to $\bdd_- A$. $\bdd$-compressing $D$ along $E$ 
      gives two disks $D_1$ and $D_2$ at least one of which has 
      boundary essential of $\bdd M$, say $D_1$. However $|D_1 \cap P|<|D \cap P|$, a 
      contradiction.
      
       Suppose that $E$ is a $\bdd$-compressing disk for $D_A$ and $E$ 
      is adjacent to $P$. Let $\alpha = E \cap D_A$. Use $E$ to isotope 
      $D$ so that a neighborhood of $\alpha$ lies in $B$, call this 
      new disk $D^1$ and let $D^1_A=D^1 \cap A$ and $D^1_B=D^1 \cap B$. 
      Note that $D^1_A$ is obtained from
      $D_A$ by cutting along $\alpha$ and $D^1_A$ is also 
      c-incompressible. Repeat the above operation naming each 
      successive disk $D^i$ until the resulting 
      surface $D^n_A= D^n \cap A$ is $\bdd$-incompressible. 
      By Corollary \ref{cor:components} $D^n_A$ consists of vertical 
      annuli and disks.  
      
      Suppose some component of $D_A$ is $\bdd$-compressible but not 
      adjacent to $\bdd_-A$. In this case the result of maximally  
      $\bdd$-compressing this component has to be a collection of 
      disks. By Case 1 of Lemma 
      \ref{lem:boundarynumber}, $|D^n_A 
      \cap P| < |D_A \cap P|$ contradicting our choice of $D$. Thus 
      every boundary compressible component of $D_A$ is adjacent to 
      $\bdd_- A$, in particular $\bdd D \subset \bdd_- A$ and $D_A$ 
      has a unique $\bdd$-compressible component $F$. By the 
      minimality assumption and Case 2 of Lemma \ref{lem:boundarynumber}, 
      some component of $D^n_A$ must be a disk. $D^n_B$ is then a 
      planar surface that we have shown must be c-incompressible and 
      has a component that is not a disk. As $\bdd D \cap \bdd_- B 
      =\emptyset$, it follows that some component of $D^n_B$ is 
      $\bdd$-compressible into $P$ and disjoint from $\bdd_- B$. 
       The above argument applied to $D^n_B$ leads to an isotopy of the 
      disk $D$ so as to reduce $D \cap P$ contrary to our assumption.
      Thus $D_A$ and $D_B$ are both collections of vertical annuli 
      and disks so $D$ is a reducing disk for $\bdd 
      M$  that intersects $P$ in a single essential 
      simple closed curve.

      \end{proof}
      
      \begin{cor}  \label{cor:stronglyirreducible}
	  Let $A\cup_P B$ be a c-strongly irreducible c-Heegaard splitting of 
	  $(M,T)$ and
	  let $F$ be a component of $\bdd M$. If $F_T$ is not parallel to $P_T$, 
	  then $F_T$ is incompressible.

\end{cor}	  
 \begin{proof}
     Suppose $D$ is a compressing disk for $F_T \subset \bdd_- A$ say. 
          By Theorem \ref{thm:one} we can take $D$ such that $|D \cap
     P|=1$,  $D_A =D \cap A$ is a compressing disk for $P_T$ lying in
     $A$ and $D-D_B$ is a vertical annulus disjoint from $T$.  As $F_T$ 
     is not parallel to $P_T$, there is a c-disk for $P_T$
     lying in $A$, $D_A$.  By a standard innermost disk and outermost arc
     arguments, we can take $D_A$ so that $D_A \cap D =\emptyset$.  But
     then $D_A$ and $D_B$ give a pair of c-weakly reducing disks for
     $P_T$ contrary to our hypothesis.
 
\end{proof}
\section{Main Theorem}

Following \cite{HS1} we will call a c-Heegaard splitting $\mcH$ {\em c-slim} 
if each component $W_{ij}=C_i \cup C_j$ obtained by cutting $M$ 
along $\mcH_-$ is c-strongly irreducible and no proper subset of 
$\mcH$ is also a multiple c-Heegaard splitting for $M$. Suppose $\mcH$ is a c-thin 
c-Heegaard splitting of $M$. If some proper subset of $\mcH$ is also a c-Heegaard 
splitting of $M$, then this c-Heegaard splitting will have lower width 
than $\mcH$. If some component $W_{ij}$ of $M-\mcH$ is c-weakly reducible, 
applying the untelescoping operation described in Lemma \ref{lem:thinning} 
to that component produces a c-Heegaard splitting of lower width. 
Thus if $\mcH$ is c-thin, then it is also c-slim.

\begin{thm} \label{thm:essentialexists}
 Suppose $M$ is a closed orientable irreducible 3-manifold containing a link $K$. If $P$ is a 
 c-strongly compressible bridge surface for $K$
  then one of the following is satisfied:

\begin{itemize}
    
    \item $P$ is stabilized
    \item $P$ is meridionally stabilized
    \item $P$ is perturbed
    \item a component of $K$ is removable
    \item $M$ contains an essential meridional surface $F$ such that $2-\chi(F_K) \leq 
2-\chi(P_K)$. 
 
 \end{itemize}   
   
 \end{thm}
\begin{proof}
Let $\mcH$ be a c-slim Heegaard splitting obtained from $P$ 
by untelescoping as in Lemma \ref{lem:thinning}, possibly in several steps. Let 
$\mcH_-$ and $\mcH_+$ denote the negative and positive surfaces of 
$\mcH$ respectively and let $W_{ij}$ be the components of $M -\mcH_-$ 
where $W_{ij}$ is the union of c-compression 
bodies $C_i$ and $C_j$ along $H_{ij}=\bdd_+C_i=\bdd_+C_j$. Suppose some component of $\mcH_-$ is 
compressible with compressing disk $D$. By taking an innermost on $D$ 
circle of $D \cap \mcH_-$ we may assume that $\bdd_- C_i$ is 
compressible in $W_{ij}$.  By Corollary \ref{cor:stronglyirreducible} 
this contradicts our assumption that $\mcH$ is c-slim. We conclude 
that $\mcH_-$ is incompressible. 

If some component of $F_K$ of $\mcH_-$ is neither an inessential 
sphere nor boundary parallel in $M_K$, then it is essential and $2-\chi(F_K) \leq 
2-\chi(P_K)$. If every component is either an inessential sphere in 
$M_K$ or boundary parallel, then by 
Lemma \ref{lem:thinning} the splitting is perturbed, stabilized, 
meridionally stabilized or there is a removable component.

\end{proof}

  \section*{Acknowledgment}
I would like to thank Martin Scharlemann for many helpful
conversations.

     \end{document}